\documentclass[12pt]{article}
\usepackage{amsmath}
\usepackage[ top=2.5cm, bottom=2.5cm, left=2.5cm, right=2.5cm ]{geometry}
\begin{document}
\linespread{1.2}
\newtheorem{theorem}{Theorem}
\newtheorem{proposition}{Proposition}
\newtheorem{corollary}{Corollary}
\newtheorem{lemma}{Lemma}

\newcommand{\al}{\alpha}
\newcommand{\si}{\sigma}

  
\title{\bf Improved Bounds on the Probability of a Union 
  and on  the Number of Events that Occur}
\author{Ilan Adler$^*$, Richard M. Karp
\thanks{Department of Industrial Engineering and Operations Research,
University of California, Berkeley
(ilan@berkeley.edu, karp@cs.berkeley.edu).
}
\and \hspace{-8.5mm} and Sheldon M. Ross
\thanks{Department of Industrial and Systems Engineering,
University of Southern California. 
 Supported  by the National Science Foundation under contract/grant   CMMI2132759. (smross@usc.edu) 
}}

 \date{April 15, 2025}

\maketitle

\begin{abstract}  
Let $A_1, A_2, \ldots, A_n$  be events in a sample space. Given the probability of the intersection of each collection of up to $k+1$ of these events, what can we say about the probability that at least $r$ of the events occur? This question dates back 
to Boole in the 19th century, and it is well known that the odd partial sums of the Inclusion- Exclusion formula provide upper bounds, while the even partial sums provide lower bounds. We give a combinatorial characterization of the error in these bounds
and use it to derive a very simple proof of the strongest possible bounds of a certain form,  as well as a couple of  improved  bounds. 
The new bounds use more information than the classical Bonferroni-type inequalities, and are often sharper. 

\end{abstract}
\section{Introduction}
Let $A_1, A_2, \ldots, A_n$  be events in a  probability space and let $X $ denote  the number of these events that occur. We are interested in obtaining bounds on $P(X \geq r)$, for $ r=1,2, \ldots, n$,   
based on knowledge of the probabilities of all  s-way intersections of the events, for each  $s=1,2, \ldots, k+1,$  for some $ k < n$.
Let $S_j$  be the sum of the probabilities of  all $j$-way intersections, for $1 \leq j \leq n.$  That is, 
$S_j = \sum_{1 \leq	 i_1 < i_2 < \cdots < i_j \leq n} P( A_{i_1} A_{i_2} \cdots A_{i_j}).$ 
 The Principle of Inclusion and Exclusion states that
\begin{equation}
P(X \geq r) = \sum_{j=r}^n (-1)^{r+j} {j-1 \choose r-1} S_j    \nonumber  
\end{equation}
It is also known that truncations of the preceding expression give upper and lower bounds on $P(X \geq r)$:  if $i$  is odd then
\begin{equation}
\label{eq:inc_exc}
 \sum_{j=r}^{r+i} (-1)^{r+j} {j-1 \choose r-1} S_j  \leq  P(X \geq r)  \leq  \sum_{j=r}^{r+i+1} (-1)^{r+j} {j-1 \choose r-1} S_j 
       \end{equation}
Many improvements of this bound have been given, and are generically called Bonferroni-type inequalities. In the next section we recall how the preceding inequalities can be derived from an algebraic identity. Our strategy
in this note is to follow a similar line of reasoning, generalizing the algebraic identity  to obtain some interesting expressions for $P(X \geq r),$ and then derive bounds from those expressions that are stronger or more general
 than the bounds previously known. Whereas our strongest bounds are expressed in terms of the individual 
$s$-fold intersections of events, and hence exploit more information than
the typical Bonferroni-type inequalities, which depend only on the coarser information given by the quantities $S_j, j = 1, \ldots, k+1,$ we also give bounds just depending on the $S_j.$  Other bounds based on the same information as ours are 
given in \cite{Gala96} and \cite{YAT2017} as well as in the references therein.\\
\vspace{.1mm}\\
{\bf{Notational Conventions.}} For integers $s,t,$ we let  ${t \choose s} = 0 $ if either $\min(s, t)< 0$ or $s>t.$ Otherwise ${t \choose s} = \frac{t!}{s! (t-s)!}.$  Also, for any event $B,$ we let $I\{B\}$ be the indicator 
variable of the event $B$, equal to $1$ if $B$ ocurs and to $0$ otherwise.

\section{Basic Expressions for the Distribution of $X$}
\begin{lemma} 
\label{lem:idendity}
For all integers $r > 0, m \geq 0$ 
\begin{equation}
 \label{eq:id1}
I\{m \geq r\} = \sum_{j=r}^m (-1)^{r+j} {j-1 \choose r-1} {m \choose j}   
       \end{equation}
       \end{lemma} \
{\bf{Proof:}}  Note that both sides of Equation $ ( \ref{eq:id1}) $ equal $0$ when $m < r.$ The proof when $m \geq r$ is by induction on $r$.  
 The result follows when $r=1$ since for $m \geq 1$  the identity 
$ \sum_{j=0}^m  {m \choose j}  (-1)^j  =  (1-1)^m = 0$  implies that 
$1 = \sum_{j=1}^m  {m \choose j}  (-1)^{j+1}  .$   
  \allowdisplaybreaks
Assuming the result for $r$, we have
     \begin{eqnarray}
 \sum_{j=r+1}^m (-1)^{r+j +1} {j-1 \choose r} {m \choose j}  &=&  \sum_{j=r+1}^m (-1)^{r+j +1} 
 [ {j \choose r}  - {j-1 \choose r-1} ]{m \choose j}  \nonumber   \\
 &=&  \sum_{j=r+1}^m (-1)^{r+j +1}  {j \choose r}  {m \choose j} + \sum_{j=r+1}^m (-1)^{r+j } 
 {j-1 \choose r-1} {m \choose j}   \nonumber  \\
  &=&  \sum_{j=r+1}^m (-1)^{r+j +1}  {m \choose r}  {m-r \choose j -r}+1 - {m \choose r} \label{eq:intermediate} \\
  &=&  {m \choose r} \sum_{k=1}^{m-r} (-1)^{k+1} {m-r \choose k} + 1 - {m \choose r}   \nonumber \\
    &=&  {m \choose r} + 1 - {m \choose r} =1  \nonumber 
   \end{eqnarray}
   where Equation (\ref{eq:intermediate}) followed from the induction hypothesis. $\quad   \rule{2mm}{2mm} $
   \begin{theorem} \
   With $X$  being the number of the events $A_1, \ldots, A_n$ that occur, 
$$P(X \geq r) =  \sum_{j=r}^n (-1)^{r+j} {j-1 \choose r-1} S_j  $$
\end{theorem} \
{\bf{Proof:}}    Because Lemma 1 holds for all $m,$ it follows that,
$$ I\{X \geq r\} = \sum_{j=r}^X (-1)^{r+j} {j-1 \choose r-1} {X \choose j}  
= \sum_{j=r}^n (-1)^{r+j} {j-1 \choose r-1} {X \choose j}  $$
where the second equality follows because $X \leq n$ and ${X \choose j} = 0 $ when $j > X.$ 
Taking expectations, and using that  $ E[{X \choose j}  ] = S_j$ yields the result. $\quad   \rule{2mm}{2mm} $
\begin{lemma}  
     \label{lem:basic_identity}
      For all nonnegative integers $m, k, r $ such that  $ k \geq r$,
$$
      \sum_{j=k+1}^m (-1)^j { j-1 \choose r-1} {m \choose j} = (-1)^{k+1} \sum_{i=1}^r   { k-i \choose r-i} {m -i \choose k-i+1 }
$$
\end{lemma} 
   {\bf{Proof:}} The result holds if $m < k$ since both sides equal $0$ in this case. Hence, we need to prove the result when $m \geq k \geq r.$ 
   The proof is by induction on $m$. As it is immediate for $m=k$ assume it is true for $m-1$. Then, with  $\;H(m,k,r) = \sum_{j=k+1}^m (-1)^j { j-1 \choose r-1} {m \choose j},$
   \allowdisplaybreaks
     \begin{eqnarray*}
     H(m,k,r) &=&  \sum_{j=k+1}^m (-1)^j { j-1 \choose r-1} [ {m -1\choose j -1} + 
     {m -1\choose j }] \\
     &=&  \sum_{j=k+1}^m (-1)^j { j-1 \choose r-1} {m -1\choose j -1} +  \sum_{j=k+1}^{m-1} (-1)^j { j-1 \choose r-1} {m-1\choose j }   \\
       &=&  \sum_{j=k }^{m-1} (-1)^{j+1} { j \choose r-1} {m -1 \choose j}  -  \sum_{j=k+1}^{m-1}
       (-1)^{j+1} { j-1 \choose r-1} {m -1\choose j }   \\
       &=&     \sum_{j=k+1 }^{m-1} (-1)^{j+1} [ { j \choose r-1}   - { j-1 \choose r-1} ]       {m -1 \choose j} \; + \; (-1)^{k+1} { k \choose r-1} {m -1 \choose k }   \\
          &=&     \sum_{j=k+1 }^{m-1} (-1)^{j+1} { j-1 \choose r-2}       {m -1 \choose j} \; + \;(-1)^{k+1} { k \choose r-1} {m -1 \choose k }   \\
                  &=&     \sum_{j=k+1 }^{m-1} (-1)^{j+1} { j-1 \choose r-2}       {m -1 \choose j} \; + \;(-1)^{k+1} [ { k-1 \choose r-1} +     { k-1 \choose r-2 } ]{m -1 \choose k }   \\
    &=&     \sum_{j=k }^{m-1} (-1)^{j+1} { j-1 \choose r-2}       {m -1 \choose j} \; + \;(-1)^{k+1}  { k-1 \choose r-1}  {m -1 \choose k }   \\
    &=&  - H( m-1, k-1,r-1)  \; + \;(-1)^{k+1}  { k-1 \choose r-1}  {m -1 \choose k }   \\
&=& (-1)^{k+1} \sum_{i=1}^{r-1}  { k-1 - i \choose r - 1 -i} {m-1 -i \choose k-1  -i+1 }  \; + \;(-1)^{k+1}  { k-1 \choose r-1}  {m -1 \choose k }   \\
&=& (-1)^{k+1} \sum_{j=2}^{r}  { k-j  \choose r - j } {m-j \choose k - j+1 }  \; + \;  ( -1)^{k+1}  { k-1 \choose r-1}  {m -1 \choose k }   \\
&=& (-1)^{k+1} \sum_{j=1}^{r}  { k-j  \choose r - j } {m-j \choose k - j+1 }  
     \end{eqnarray*}
     where the induction hypothesis was used to evaluate $H(m-1, k-1, r-1).  \quad   \rule{2mm}{2mm} $
    \begin{theorem}
  \label{thm:bounds_k}
 For $k \geq r,$
$$ P(X \geq r) =  \sum_{j=r}^k (-1)^{r+j} {j-1 \choose r-1} S_j  + (-1)^{r+k+1}  \sum_{i=1}^r {k-i \choose r-i}  E[{X-i \choose k-i+1}]    $$
     {\bf{Proof:}} It follows from   Lemmas 1 and 2 that for $k \geq r$
     $$ I\{X \geq r\} = \sum_{j=r}^k (-1)^{r+j} {j-1 \choose r-1} {X \choose j}   +  (-1)^{r+k+1} \sum_{i=1}^r   { k-i \choose r-i} {X -i \choose k-i+1 } $$
     Taking expectations now yields the result. $  \quad   \rule{2mm}{2mm} $
         \end{theorem}
     {\bf{Remark:}} The classical Bonferroni type inequalities \eqref{eq:inc_exc} follow from Theorem 2 upon observing  that
$ \, \sum_{i=1}^r {k-i \choose r-i} E[{X-i \choose k-i+1}]   \geq 0. $
  \section{Bounds on $P(X \geq r)$}
    \begin{lemma}
  \label{lem:BB}
   $\;E[{ X-i \choose  k+1-i} ] \geq     \frac{  { k+1 \choose  i} }  {  {n \choose  i} } S_{k+1} $
    \end{lemma}
  {\bf{Proof:}}    Because $X \leq n$ and  $ { j-i \choose  k+1-i} / { j \choose  k+1}  $ is monotonically decreasing in $j$,  it follows that 
$${ X-i \choose  k+1-i} \geq   { X \choose  k+1}   \frac{  { n-i \choose  k+1-i} } {  { n \choose  k+1} }= 
{ X \choose  k+1}  \frac{  { k+1 \choose  i} } {  {n \choose  i} } $$
 Taking expectations yields the result.  $  \quad  \rule{2mm}{2mm} $ 
 \begin{theorem}  
\label{thm:Benf}
$$ P(X \geq r) \geq  \sum_{j=r}^k (-1)^{r+j} {j-1 \choose r-1} S_j  +   \sum_{i=1}^r {k-i \choose r-i}   \frac{  {k+1 \choose i}} {  {n \choose i} } S_{k+1} ,\quad  r+k \;\mbox{ odd} $$
  $$ P(X \geq r) \leq  \sum_{j=r}^k (-1)^{r+j} {j-1 \choose r-1} S_j  -   \sum_{i=1}^r {k-i \choose r-i}   \frac{  {k+1 \choose i}} {  {n \choose i} } S_{k+1} ,\quad  r+k \;\mbox{ even} $$
   \end{theorem}
  {\bf{Proof:}} Immediate from Theorem \ref{thm:bounds_k} and Lemma \ref{lem:BB}.  $  \rule{2mm}{2mm} $ \\
  \vspace{.1mm}\\
  {\bf Remarks:}
\begin{itemize}
\item
Note (by considering the case $P(X=n)=1$) that $\alpha=\sum_{i=1}^r {k-i \choose r-i}  \frac{  {k+1 \choose i}} {  {n \choose i} } $ is the most stringent among all the bounds of the type
$$ \sum_{j=r}^k (-1)^{r+j} {j-1 \choose r-1} S_j  + (-1)^{r+k+1}\alpha S_{k+1} $$
\item
In his book, Galambos \cite{Gala96} presents a  result (proved first in \cite{Margaritescu87}) that is similar to Theorem 3, except that  $\sum_{i=1}^r {k-i \choose r-i}  \frac{  {k+1 \choose i}} {  {n \choose i} } $ is
 replaced by $\frac{( \sum_{j=0}^{k-r} (-1)^j {r+j-1 \choose r-1} {n \choose r+j}-1)}{{n \choose k+1}}$. However, these terms are actually identical, as can be easily shown via Lemmas \ref{lem:idendity}, \ref{lem:basic_identity},  and  \ref{lem:BB}. 
 The proof presented here, though, is much simpler.
\end{itemize}
Considering Theorem \ref{thm:bounds_k}, it is clear that obtaining a lower bound for  
$
E\left[\binom{X-i}{k+1-i} \right]
$ 
can provide upper and lower bounds for \( P(X \geq r) \). 
So far, we have derived bounds based solely on the coarse information provided by the quantities  $S_j, j = 1, \ldots, k+1,$ 
 Next, we refine these bounds by expressing them in terms of the individual  s-fold intersections of events, leading to stronger bounds.

\medskip

  Let  ${\cal{S}} =\{A_{i_1}, \ldots, A_{i_x} \}, i_1 < \ldots < i_x, $ be the set of events that occur,
  and note that for $i  \leq k \leq x-1$, 
  $ { x-i \choose  k+1-i} $ is the number of $k+1-i$ sized subsets in . 
  $\{A_{i_1}, \ldots, A_{i_{x}-i}\}$. 
Now, for $r_1 < \ldots <  r_{k+1-i}$,  $\{A_{r_1}, \ldots, A_{ r_{k+1-i}}\}$ is such a subset  if all of the events
$A_{r_1},  \ldots, A_{ r_{k+1-i}}$ occur as does at least $i$ events $A_{t_1}, \ldots, A_{t_i}$, where $r_{k+1-i }< t_j  \text{ for }j=1, \ldots, i$.
Hence, defining  
\[
T(  r_1 , \ldots , r_{k+1-i}) = \{(t_1, \ldots, t_i):  r_{k+1-i} < t_1 < \cdots <  t_i \}
\]
we obtain the expectation bound  
\begin{align}
E\left[\binom{X - i}{k+1-i}\right] &=  
\sum_{\substack{r_1 < \cdots < r_{k+1-i}}}  
P\left( \bigcup_{(t_1, \ldots, t_i) \in T(r_1, \ldots, r_{k+1-i})}  
A_{r_1}  \cdots A_{r_{k+1-i}} A_{t_1} \cdots A_{t_i} \right) \nonumber \\
 \label{eq:r_max}
&\geq  
\sum_{\substack{r_1 < \cdots < r_{k+1-i}}}  
\max_{(t_1, \ldots, t_i) \in T(r_1, \ldots, r_{k+1-i})}  
P(A_{r_1} \cdots A_{r_{k+1-i}} A_{t_1} \cdots A_{t_i}).\end{align}

\noindent The preceding yields
\footnotesize
 \begin{theorem}   
\label{thm:r_max}
\begin{equation*}
{P(X \geq r) \geq  \sum_{j=r}^k (-1)^{r+j} {j-1 \choose r-1} S_j  +  
 \sum_{i=1}^r {k-i \choose r-i}   \sum_{r_1  < \cdots < r_{k+1-i}}  \max_{(t_1 \ldots, t_i )\in T(r_1, \ldots, r_{k+1-i})} P(A_{r_1} \cdots  A_{r_{k+1-i }} A_{t_1} \cdots  A_{t_i } )}
 \;r+k \mbox{ odd} 
\end{equation*}
\begin{equation*}
{P(X \geq r) \leq  \sum_{j=r}^k (-1)^{r+j} {j-1 \choose r-1} S_j  -  
 \sum_{i=1}^r {k-i \choose r-i}   \sum_{r_1  < \cdots < r_{k+1-i}}  \max_{(t_1 \ldots, t_i )\in T(r_1, \ldots, r_{k+1-i})} P(A_{r_1} \cdots  A_{r_{k+1-i }} A_{t_1} \cdots  A_{t_i } )}
 \;r+k \mbox{ even} 
\end{equation*}
   \end{theorem}
   \normalsize
  {\bf{Proof:}} Immediate fromTheorem \ref{thm:bounds_k} and \eqref{eq:r_max}$. \quad   \rule{2mm}{2mm} $ 
  \subsection{An Improved 
  Bound of $P(X \geq 1) $}
 
The case of $P(X \geq 1) $ is of particular interest, as there are applications in which the focus is on the probability of at least one event occurring (see e.g. \cite{Gala96}).
 Specializing $i=1$ in \eqref{eq:r_max}, we obtain
\begin{align}
E\left[\binom{X - 1}{k}\right] 
\geq  
\sum_{r_1  < \cdots < r_{k}}  \max_{r>r_k } P(A_{r_1} \cdots  A_{r_{k }} A_{r}  )
\label{eq:r1_permu}
\end{align}

As the preceding lower bound holds for all numberings of the $n$ events,  an  improved lower bound can be obtained by finding the numbering that maximizes the right-hand side of \eqref{eq:r1_permu}. However, this is not computationally feasible when 
n is large. Instead, we will show how to compute the expected value of the right-hand side of \eqref{eq:r1_permu} when the numbering is chosen at random. To this end, we first present the following lemma.
      \begin{lemma}
     \label{lem:perm} 
Let $i_1, \ldots, i_k, j_1, \ldots, j_s$ all be distinct indices. The probability that $j_s$ is, in a random permutation, the only one of $j_1, \ldots, j_s$ to follow all of the indices $i_1, \ldots, i_k$ is  $\; \frac{k}{(k+s)(k+s-1)}.$ 
 \end{lemma}
{\bf{Proof.}}  Note that this will be true if and only if $j_s$ is the last and one of $i_1, \ldots, i_k$ is the next to last of these $k+s$ indices to appear in the random permutation.  Consequently, the  desired probability is  $\;\frac{1}{k+s}  \frac{k}{k+s -1} .  \quad   \rule{2mm}{2mm} $ \\
     \vspace{.1mm}\\
For distinct indices $i_1, \ldots, i_k$, let $W_s(i_1, \ldots, i_k)$ be the $s^{th}$ largest of the $n-k$ values 

\noindent$P(A_{i_1} \cdots A_{i_k} A_j),\, j \notin \{i_1, \ldots, i_k\}.$ 
    \begin{lemma}
    \label{lem:w_bound}
      $$ E[ {X - 1 \choose k}] \,  \geq  \sum_{i_1 < \cdots < i_{k}}  \sum_{s=1}^{n-k} W_s(i_1, \ldots, i_k) \frac{k}{(k+s)(k+s-1)}
  $$
  \end{lemma}
     {\bf{Proof:}}  
With  $p_1, \ldots, p_n$ being  a random permutation of $1, \ldots, n,$ and  $B_i = A_{p_i}, i = 1, \ldots, n,$ it follows from $ (\ref{eq:r_max})$ that 
\begin{equation}  \label{r1}
 E[ {X - 1 \choose k}]  \geq  \sum_{i_1 < \ldots < i_k} E[  \max_{r>i_k} P(  B_{i_1} \cdots B_{i_k}B_r) ] 
\end{equation}
Now, for $i_1 < \cdots < i_k$,
\begin{equation}  \label{r2}
 E[  \max_{r>i_k} P(  B_{i_1} \cdots B_{i_k}B_r) ] = \frac{1}{{n \choose k}} \,\sum_{r_1 < \cdots < r_k} E[\max_r P(A_{r_1} \cdots A_{r_k} A_r) ]
 \end{equation}
 
\noindent where the maximum on the right-hand side of the preceding equation (as well as on the left-hand side of the succeeding  equation) is taken over all indices $r$ that follow all of $r_1, \ldots, r_k$ in the random permutation.  Applying  Lemma \ref{lem:perm}  gives that 
\begin{equation}  \label{r3}  
E[\max_r P(A_{r_1} \cdots A_{r_k} A_r) ] = \sum_{s=1}^{n-k} W_s(r_1, \ldots, r_k) \frac{k}{(k+s)(k+s-1)}
 \end{equation}
 The result now follows from $(\ref{r1}), (\ref{r2}), $ and $(\ref{r3}).\quad   \rule{2mm}{2mm} $ 
 
 \medskip
 
\noindent The preceding yields 
\begin{theorem}
  \label{thm:av_perm}
  $$P(X \geq 1) \geq \sum_{j=1}^{k} (-1)^{j+1} S_j  +  \sum_{r_1 < \cdots < r_{k}}  \sum_{s=1}^{n-k} W_s(r_1, \ldots, r_k) \frac{k}{(k+s)(k+s-1)} \quad \mbox{$\;k\;$  even} $$
    $$P(X \geq 1) \leq \sum_{j=1}^{k} (-1)^{j+1} S_j - \sum_{r_1 < \cdots < r_{k}}  \sum_{s=1}^{n-k} W_s(r_1, \ldots, r_k) \frac{k}{(k+s)(k+s-1)} \quad \mbox{$\;k\;$ odd} $$
    \end{theorem}
     {\bf{Proof:}} Immediate fromTheorem \ref{thm:bounds_k} and Lemma \ref{lem:w_bound}. $  \quad   \rule{2mm}{2mm} $ 

\bigskip 

      {\bf{Example}} Suppose that 
        $n=6$ and that $P(A_{i_1} \cdots A_{i_j})=   \al_{i_1} \cdots \al_{i_j}, j \leq k+1$, where 
     $\al_t = \frac{t+18}{100}, t=1, \ldots, 6$, so $S_1 = 1.290 , S_2 =   0.6925   , S_3 = 0.1980   $. 
     
     Below, we tabulate the various bounds  presented in this paper for $r=1$ and $k=2$.  
     
     \noindent Note that since $k$ is even ,   $$S_1 - S_2 + (\text{lower bound on } E[ {X -1 \choose 2}]) \leq P(X  \geq 1) \leq   S_1 - S_2 + S3$$

 \medskip
 
\scriptsize 
\begin{center}
        \begin{tabular}{||c||c||c|c|cc||}
\hline
\hline
&&&
  & &\\
Reference  &$S_1 - S_2$  & Lower Bound &
Lower Bound  & Upper Bound&\\
& & on $E\left[\binom{X - 1}{2}\right]$&
 on $P(X \geq 1)$ & on $P(X \geq 1)$&\\
   && & &($S_1-S_2 + S_3$) & \\  
\hline
\hline
& & & & &\\
Eq. \eqref{eq:inc_exc} (Bonferroni)& .5975 & 0 & .5975 &.7955 &\\
  && & & &\\ 
\hline
& & & & &\\
  Theorem \ref{thm:Benf}    &   .5975 &
 .0990  &  .6965  &.7955 &\\
  && & & &\\ 
 \hline
 & & & & &\\ 
  Theorem \ref{thm:r_max}&  .5975   &
.1057  &.7032 &.7955 &\\
& & & & &\\ 
\hline
   & & & & &\\ 
    Theorem \ref{thm:av_perm}   &  .5975  &
.1896  & .7871& .7955&\\
    & & & & &\\ 
\hline
\hline
\end{tabular}
\end{center}


\bigskip

\normalsize

\medskip

 \begin{thebibliography}{1}
 
 \bibitem{Gala96}
J. Galambos and I. Simonelli.
\newblock {\em Bonferroni-type Inequalities with Applications}.
\newblock Springer-Verlag,  New York, 1996.

\bibitem{Margaritescu87}
E. Margaritescu
\newblock {\em On some Bonferroni inequalities}.
\newblock {  Stud. Cerc. Mat., 39 246-251, 1994}.

\bibitem{YAT2017}
J. Yang, F. Alajaji and G. Takahara.
\newblock {\em A short survey on bounding the union probability using partial information}.
\newblock { arXiv:1710.07576, 2017}.
\end{thebibliography}
      \end{document}